\def\ZZ         {{\bf Z}}

\def\CC         {{\bf C}}
\def\QQ         {{\bf Q}}
\def\PP         {{\bf P}}

\def\L          {{\cal L}}
\def\A          {{\cal A}}
\def\S          {{\cal S}}
\def\O         {{\cal O}}
\def\F          {{\cal F}}

\documentstyle[twoside,12pt]{article}
\newtheorem{prop}{Proposition}[section]

\newtheorem{theo}[prop]{Theorem}

\newtheorem{rem}[prop]{Remark}
\newtheorem{coro}[prop]{Corollary}

\newtheorem{exam}[prop]{Example}

\title {\bf Eigenvalues for the monodromy of the Milnor fibers of
arrangements}
\author {
{\sc Anatoly Libgober}\\
{\small\it Department of Mathematics, University of Illinois,
Chicago, Ill 60607}}

\begin{document}

\date{\today}
\maketitle
\begin{abstract}
{We decribe upper bounds for the orders of the eigenvalues of the 
monodromy of Milnor fibers of arrangements given in terms 
of combinatorics.} 
\end{abstract}

\section{\bf Introduction.}

The central object is the study of the topology of isolated
hypersurface singularities is the Milnor fiber. If
$f(x_0,...,x_n)$ has an isolated singularity, say at the origin
$\cal O$, then the Milnor $F_f$ of $f$ is the intersection of
$f=t$ with a ball $B_{\epsilon}$ of a small radius $\epsilon$
about the origin ($\vert t \vert \ll \epsilon$). The Milnor fiber
$F_f$ is an $n$-connected smooth manifold with boundary having the
natural action of the monodromy obtained by letting $t$ vary
around a small circle around the origin of $\bf C$. This is a
diffeomorphism of $F_f$ constant on the boundary which isotopy
class modulo boundary depends only on $f$. In particular, one has
a well defined operator $T_f$ on $H_n(F,{\bf Z})$. This operator
carries much of information about the topology of $f$ (cf.
\cite{milnor}) and its calculation can be done in several ways.
Particularly effective is the method due to A'Campo based on a
resolution of the singularity of $f$ (cf. \cite{acampo}).

If the singularity is not isolated the Milnor fiber may ``lose''
connectivity. The monodromy diffeomorphism is still well defined
up to isotopy and induces a well defined operator $T_{f,i}$ on
each homology group $H_i(F_f,{\bf Z})$. Though the zeta function:
$\prod_i {\rm det} (T_{f,i}-I)^{(-1)^{i-1}}$ still can be easily
calculated by A'Campo method, the eigenvalues of $T_{f,i}$ for
each $i$ reflect more subtle properties of the singularity not
necessarily encoded into the topological data about resolutions
(cf. \cite{steenbrink} \cite{withmihai}).

It is well known that the Milnor fiber of a cone over a projective
hypersurface of degree $d$ can be identified with the $d$-fold
cyclic cover of the complement to the hypersurface (cf. for
example \cite{annals} or section \ref{prelim} below). In
particular information about the homology of cyclic covers is
equivalent to the information about the homology of Milnor fibers.
The former were studied extensively over the last 20 years (cf.
\cite{duke}, the survey \cite{china} or \cite{dimca}). In
particular it was shown that the eigenvalues of the deck
transformations  acting on the homology of cyclic covers and hence
the eigenvalues of the monodromy of related Milnor fibers depend
on position of singularities and bounded by the local type of
hypersurfaces and their behavior at infinity.

In this note we shall specialize the results surveyed in
\cite{arcata} and \cite{china} to the case of arrangements. These
results fall into two groups. The first contains restrictions on
the orders of the eigenvalues (note that eigenvalues of the
monodromy are roots of unity as a consequence of the monodromy
theorem, cf. \cite{milnor}). These restrictions imposed entirely
by the combinatorics of arrangement (cf. sect. 2.4 and theorem
\ref{main}. In particular one obtains vanishing of the cohomology
of certain local systems (i.e. those which cohomology are the
components of the cohomology of cyclic covers). Our approach does
not depend on Deligne's (\cite{diffeq}) results and gives
vanishing of cohomology in specific dimensions i.e. we obtain
conditions for ''non resonance'' in certain range.

The second group of results deals with a  calculation of
characteristic polynomial of the monodromy in the case of line
arrangements in terms of dimensions of system of curves determined
by collection of vertices of the arrangement. Methods used below
are all contained in our previous work but they give more precise
results in the case when the hypersurface is an arrangement. This
case of arrangements was considered recently in (cf.
\cite{cohenorlik}, \cite{dedham}).

This work was supported by NSF grants DMS-9872025 and DMS 9803623.

\section{\bf Preliminaries.}
\label{prelim}

\paragraph{2.1.}
Let $\A$ be an arrangement of hyperplanes in $\PP^n$ and
$l_i(x_1,...,x_{n+1})=0, i=1,..,d$ be the defining equations for the
hyperplanes of $\A$. To $\A$ corresponds the cone $A$ given by the
equation  $\prod l_i=0$ in $\CC^{n+1}$. Since the defining
equation of $A$ is homogeneous it follows that the Milnor fiber
$F_A$ can be identified with the affine hypersurface $\prod
l_i=1$.

\paragraph{2.2.}
Recall that  ${H}_1(\PP^n-\A,\ZZ)=\ZZ^{d-1}$ with generators given
by the loops $\lambda_i$ each being the boundary of a small 2-disk
in $\PP^n$ transversal to $l_i=0$ (i.e. a meridian) at a
non-singular point of $\A$. These loops satisfy single relation
$\sum_{i=1}^{i=d} \lambda_i=0$. In particular $\lambda_i
\rightarrow 1 \ {\rm mod} \ d$ defined the homomorphism
${H}_1(\PP^n-\A,\ZZ) \rightarrow \ZZ / d \ZZ$. The corresponding
branched covering can be realized as the hypersurface $\bar
V_{\cal A}$ in $\PP^{n+1}$ given by
\begin{equation}
\label{milnoreq} x_0^d=\prod l_i(x_1,...,x_{n+1})
\end{equation}
The covering map is induced by the projection
\begin{equation}
\label{projection}
 \pi: \ \ \ (x_0,x_1,...,x_{n+1}) \rightarrow (x_1,...,x_{n+1}).
 \end{equation}
The hypersurface $\bar V_{\A}$ is singular with the singular locus
$Sing$ consisting of the points taken by this projection into
points in $\PP^n$ belonging to at least two hyperplanes of the
arrangement. The unbranched covering $V_{\A}$ of $\PP^n-\A$
corresponding to $H_1(\PP^n-\A,\ZZ) \rightarrow \ZZ/d\ZZ$ is the
restriction of this projection to $\bar V_{\A}-\{x_0=0 \}$. The
natural identification of $\CC^{n+1}$ with the complement to
$x_0=0$ in $\PP^{n+1}$ yields the identification of the Milnor
fiber $F_A$ with unbranched covering $V_{\A}$. Moreover the
description of the monodromy of Milnor fiber for weighted
homogeneous hypersurfaces in \cite{milnor} yields that a monodromy
diffeomorphism of $F_A$  can be identified with the deck
transformation $x_0 \rightarrow \mu_d \cdot x_0$ ($\mu_d$ is a
primitive root of unity of degree $d$) of $V_{\A}$. The operator
induced on the homology of either of these spaces we denote as
$T_A$.

An immediate consequence of this description of the monodromy is
that for an arrangement of $d$ hyperplanes each eigenvalue of
$T_{A,i}$ acting on $H_i(F_{A},\ZZ)=H_i(V_{\cal A},\ZZ)$ has an
order dividing $d$. The multiplicity of the eigenvalue $1$ for
$T_{A,i}$ acting on $H_i(F_A,\ZZ)$ is equal to ${\rm k}
H_i(\PP^n-\A,\ZZ)$. Indeed, in the cohomology spectral sequence
$H^q(\ZZ / d \ZZ,H^p(V_{\A})) \Rightarrow H^{p+q}(\PP -\A)$
associated with action of $\ZZ / d \ZZ$ all terms with $p \ge 1$
are zeros and the rank of invariant part of $H^p(V_{\A})$
coincides with the multiplicity of the eigenvalue 1 of the
monodromy.

\paragraph{2.3.} Milnor fiber of a central arrangement is homotopy
equivalent to the infinite cyclic cover
\begin{equation} \label{infinitcover}
\widetilde {({\CC}^{n+1}-A)}_{\infty} \rightarrow {\CC}^{n+1}-A
\end{equation}
corresponding to the homomorphism $\ZZ^d
=H_1({\CC}^{n+1}-A)\rightarrow \ZZ$ sending a generator
corresponding to a hyperplane to the positive generator of $\ZZ$.
Indeed, as a loop around the hyperplane $l_i=0$ one can take the
intersection of $n$ hypersurfaces
$l_1=r_1,..,l_{i-1}=r_{i-1},l_{i+1}=r_{i+1},..,l_{n}=r_n,\prod_{i
\ge n} l_i=r_{n+1}$ (equation $l_i$ omitted) and $\vert l_i \vert
= \epsilon$. This loop is taken by the map $f_A: (x_1,..,x_{n+1})
\rightarrow \prod l_i(x_1,...,x_{n+1})$ into a small circle about
the origin of $\CC^*$. Hence $\widetilde{(\CC^{n+1}-A)_{\infty}}$
is homotopic to the fiber product $\CC^{n+1}-A \times_{\CC^*} \CC$
with respect to $f_A$ and $exp: \CC \rightarrow \CC^*$. This fiber
product clearly is homotopy equivalent to $f_A^{-1}(1)$.
 Alternatively, preimage of $\CC_a \ (a \in \PP^n-\A)$ in
$\widetilde {(\CC^{n+1}-A)_{\infty}}$ of $\CC^*_a$ which a fiber
of $\CC^{n+1}-A \rightarrow \PP^n-\A$ consists of $d$ contractible
component since the image of $\pi_1(\CC^*)$ in the Galois group of
the cover (\ref{infinitcover}) is an infinite subgroup of index
$d$.

If $\widetilde{({\CC}^{n+1}-A)_e}$ is a cyclic cover corresponding
to the map $H_1({\CC}^{n+1}-A)\rightarrow \ZZ/e\ZZ$ then the
eigenvalues of the deck transformation on
$\widetilde{({\CC}^{n+1}-A)_e}$ have order which divides the least
common multiple of $d$ and $e$. Indeed, denoting for a CW-complex
 $X$ and its cyclic cover $\tilde X$ by $C_i(X)$ and $C_i(\tilde X)$
 the group of
 $i$-chains we have the exact sequence:
 $0\rightarrow
C_i(\tilde X) \rightarrow C_i(\tilde X) \rightarrow C_i(X)
\rightarrow 0$. Here the left homomorphism is the deck
transformation minus identity. This yields homology sequence (cf.
\cite{milnorcovers}):
$$
...\rightarrow H_i(\widetilde {({\CC}^{n+1}-A)}_{\infty})
\rightarrow H_i(\widetilde {({\CC}^{n+1}-A)}_{\infty}) \rightarrow
H_i(\widetilde {({\CC}^{n+1}-A)}_{e}) \rightarrow$$
\begin{equation}
\label{milnorsequence}
H_{i-1}(\widetilde {({\CC}^{n+1}-A)}_{\infty}) \rightarrow ...
\end{equation}

The $\CC[t,t^{-1}]$-module
 $H_i(\widetilde {({\CC}^{n+1}-A)}_{\infty}) $ with $t$ acting
as the deck transformation is annihilated by $t^d-1$ as a
consequence of interpretation of $\widetilde
{({\CC}^{n+1}-A)}_{\infty}$ as the Milnor fiber. On the other hand
the left map in (\ref{milnorsequence}) is multiplication by
$t^e-1$. Hence the claim follows.

As a corollary we obtain the following:

\begin{prop} \label{divisibilityofd}
Let $\widetilde {({\CC}^{n+1}-A)}_e$ be a $e$-fold
cover of $\CC^{n+1}-A$ and $\rm d$ the number of hyperplanes in
$A$. Then an eigenvalue of a deck transformation acting on
$H_i(\widetilde {({\CC}^{n+1}-A)}_e)$ has an order dividing both
$e$ and $d$.
\end{prop}

\paragraph{2.4.} In the case $n=2$ one can relate the homology of
$V_{\A}$ to the
homology of a non-singular compactification $\tilde V_{\A}$ of
$V_{\A}$ (cf. \cite{duke},\cite{loeser}):
\begin{equation} \label{firsthomology}
{\rm rk} H_1(\tilde V_{\A},\QQ)={\rm rk} H_1(V_{\A},\QQ)-(d-1)
\end{equation}
Indeed, the argument in \cite{duke} shows that ${\rm rk}
H_1(\tilde V_{\A})={rk}H_1(\bar V_{\A}-{\rm Sing})$. The rest
 follows from the exact sequence of the pair $(\bar V_{\A}-{\rm
Sing },V_{\A})$. Indeed, from the excision and Lefschetz duality one
obtains the isomorphism of
 $H_i(\bar V_{\A}-{\rm Sing },V_{\A})$ and $H^{4-i}$ of the disjoint
union of lines of $\A$ and hence
 ${H_1}(\bar V_{\A}-{\rm Sing },V_{\A})=0,
 {H_2}(\bar V_{\A}-{\rm Sing },V_{\A},\ZZ)=\ZZ^d$. Moreover,
the deck transformations act trivially on the latter group and $\ZZ/ d
\ZZ$-invariant subgroup of ${H}_2(\bar V_{\A}-Sing)$ is cyclic
which injects into ${H}_2(\bar V_{\A}-{\rm Sing },V_{\A})$.
This yields (\ref{firsthomology}).

For a space $V$ acted upon by $T$ let $V_{\xi}$ be the
subspace spanned by the eigenvectors with eigenvalue $\xi$.
The above arguments, since $H_i(F_A)=H_i(V_{\A})$, also show that
\begin{equation} \label{nonones}
\oplus_{\xi \ne 1} H_1(F_A)_{\xi}=
\oplus_{\xi \ne 1} H_1(\tilde V_{\cal A})_{\xi}, \ \ rkH_1(F_A)_1=d-1
\end{equation}

\paragraph{2.5}
Finally, recall two results which will be used below. The first is
the Lefschetz hyperplane section theorem (cf. \cite{goreski}). Let
$X$ be a stratified complex algebraic variety having the dimension
$n$ and let $H$ be a hyperplane transversal to all strata of $X$.
Then the homomorphism $H_i(X \cap H, {\bf Z}) \rightarrow
H_i(X,\ZZ)$ induced by injection $X \cap H \rightarrow X$ is an
isomorphism for $i \le n-2$ and is surjective for $i=n-1$.

Secondly, we shall use the Leray spectral sequence corresponding
to a covering. More precisely, let $X=\bigcup_{i=1}^{i=N} U_i$ be
a union of locally closed  subsets. Then there is the
Mayer-Vietoris spectral sequence:
$$E_2^{p,q}=\oplus_{i_1 <...<i_q} H^p(U_{i_1} \cap ... \cap U_{i_q})
\Rightarrow H^{p+q}(X)$$ Moreover, if a group $G$ acts on $X$
leaving each $U_i$ invariant then all the maps in this spectral
sequence are equivariant (cf. \cite{godement}).

\section{\bf Bounds on the orders of the eigenvalues.}
Let $\A$ be an arrangement in ${\bf P}^n$. Let us call two points
$P$ and $P'$ in ${\bf P}^n$ equivalent if collections of
hyperplanes in the arrangement containing $P$ and $P'$ coincide.
Each equivalence class is a smooth submanifold of $\PP^n$ and
equivalence classes form a stratification of ${\bf P}^n$ with the
union of strata of dimension $n-1$ coinciding with the union of
hyperplanes in $\A$. Let $S_1^k,...,S_{s_{k}}^k$ be the collection
of strata of codimension $k$. For each stratum we define the
multiplicity $m(S^k_i)$ as the number of hyperplanes in
arrangement containing a point from this stratum. Any hyperplane
$H \in \A$ acquires the induced stratification.

\begin{theo}
\label{main} Let $H$ be a hyperplane of the arrangement $\A$ and
let $m_1^k(H)=m(S_1^k),...,m^k_{s(k)}=m(S^k_{s(k)})$ be the
collection of multiplicities of strata of the above stratification
of $\PP^n$ which belong to $H$ and have codimension $k$ . Let
$\xi$ be an eigenvalue of the monodromy of the Milnor fiber $F_A$
of $\cal A$ acting on $H_{k-1}(F_A,{\bf C})$. Then for any $H \in
\A$ one has $\xi^{m_i^j(H)}=1$ for at least one of the
multiplicities $m^j_i(H)$ with $j \le k$.
\end{theo}

\par {\bf Proof.} The idea is to bound the orders (of the
eigenvalues of the deck transformations acting on the homology of
the cyclic covering induced by $\pi$ (cf (\ref{projection})) on
{\it a tubular neighborhood of $H$} and then use a Lefschtez type
argument to derive the theorem.
\par First, however, let us notice that it is enough to show the theorem in
the case $k=n$. Indeed if $k<n$ then for a generic linear subspace
$L \in \PP^n$ having dimension $k$ and transversal to all strata
of $\A$ we have the induced arrangement $\L=\A \cap L$ which has,
due to transversality, as the multiplicities of its strata the
integers $m_i^j, j \le k, 1 \le i \le s_j $. Moreover, one has an
equivariant with respect to the group of deck transformations map:
$H_j(V_{L \cap \A}) \rightarrow H_j(V_{\A})$, which, by Lefschetz
theorem, is an isomorphism for $j \le k-2$ and surjective for
$j=k-1$. Hence the above theorem for the arrangement $L \cap \A$
yields the result in general and we shall assume from now on that
$k=n$.
\par  Let $T(H) \subset \PP^n$ be a small tubular
neighborhood of $H$, $B=H-H \cap \A$ and $T(B)=T(H)-\A$. The above
stratification of $\PP^n$  yields a partition of $T(B)$ into union
of subsets of ${\bf P}^n$ corresponding to the strata of the above
stratification of $H$ so that each subset is a locally trivial
fibration over the corresponding stratum in $H$. The fiber of this
fibration over a stratum of dimension $k$ is a central arrangement
of hyperplanes in ${\bf C}^{n-k}$ and the number of hyperplanes in
the latter is equal to the multiplicity of the stratum. Let ${\cal
S}^k_l$ be a collection of subsets of $T(B)$ each of which fibers
over the stratum $S^k_l$ and chosen so that their union in $T(B)$.

Intersection of subsets $\S^{k_1}_{l_1}$ and $\S^{k_2}_{l_2}$ such
that $k_1 \ge k_2$ is non-empty if and only if the stratum
$S^{k_2}_{l_2}$ is in the closure of $S^{k_1}_{l_1}$. In this case
this the intersection is the fibration with the same fiber as
$\S^{k_1}_{l_1} \rightarrow S^{k_1}_{l_1}$ and the base being a
subset in the stratum $S^{k_1}_{l_1}$ (more precisely the base is
the complement in a small neighborhood of the closure of
$S^{k_2}_{l_2}$ in the closure of $S^{k_1}_{l_1}$ to the union of
the hyperplanes of the arrangement induced by $\cal A$ on the
closure of $S^{k_1}_{l_1}$)
. \footnote{To illustrate this, let us
consider $\A$ which is the set of zeros of $x\cdot y \cdot z \cdot
$ in $\PP^3$. We have one stratum of dimension 3, four strata of
dimension 2 corresponding to four planes of the arrangement, six
one dimension strata and corresponding to lines and strata of
dimension zero. Sets $\S^2_i$ are fibrations over $\PP^2$ minus
three lines and having as a fiber the circle. Sets $\S^1_i$ are
fibered with the fiber homeomorphic to $\CC^2$ minus a pair of
intersecting lines. The base is $\CC$ minus a point, etc.
Intersection of strata having as their closure a plane and a line
in this plane is the fibration with the fiber being a circle over
a regular neighborhood of $\PP^1$ in $\PP^2$ minus two of its
fibers.}

We shall denote the $\pi$-preimage of each of the sets $\S_i^k$
(resp. $T(B)$) as $\tilde \S_i^k$ (resp. $\widetilde {T(B)}$).

\par Let us consider the Mayer-Vietoris
spectral sequence (cf. 2.5):
\begin{equation}
\label{mayervietoris} E_2^{p,q}: \oplus H^p(\tilde \S^{k_1}_{t_1}
\cap ...\cap \tilde \S^{k_{q+1}}_{t_{q+1}}) \Rightarrow
H^{p+q}(\widetilde {T(B)})
\end{equation}

The sequence (\ref{mayervietoris}) is equivariant with respect to
the action of the group of deck transformations of the cover $\pi$
restricted to $\widetilde {T(B)}$. An eigenvalue of the deck
transformation acting on $E_2^{p,q}$ must satisfy:
$\xi^{m^i_{t_i}}$ where $i \le k$. To see  this notice that each
summand in (\ref{mayervietoris})  fibers over a subset of $S^i_j$
with the fiber being the cyclic cover of an arrangement of $m^i_j$
hyperplanes. Consider the Leray spectral sequence for such
fibration:
\begin{equation} \label{leray}
H^a(S^{k_1}_{t_1} \cap ...\cap
S^{k_l}_{t_l},H^b(\widetilde{\F(\S^{k_1}_{t_1} \cap ...\cap
\S^{k_l}_{t_l})})) \Rightarrow H^{a+b}(\tilde \S^{k_1}_{t_1} \cap
...\cap \tilde \S^{k_l}_{t_l})
\end{equation}
 where $\widetilde{\F(S^{k_1}_{t_1}
\cap ...\cap \S^{k_l}_{t_l})}$ is the cover of the arrangement
$\F(S^{k_1}_{t_1} \cap ...\cap \S^{k_l}_{t_l})$ which is the fiber
of the fibration associated with corresponding stratum and hence
is an arrangement of $m^i_{t_i}$ hyperplanes where $i \le p$.
(\ref{leray}) yields that the degree of an eigenvalue on
$H^p(S^{k_1}_{t_1} \cap ...\cap \S^{k_l}_{t_l})$ is $m^i_{t_i}$
and the claim follows.
\par Now let us consider a generic (i.e. transversal to all strata)
hypersurface
${\cal V}$ which belongs to $T(B)$: hypersurface
$x_0^N=\tilde B$ where $\tilde B$ is a sufficiently small
deformation of $B$ is a possible choice. Since ${\cal V}$ is an
ample divisor in ${\bf P}^{n+1}$ by Lefschetz theorem (\cite{nori})
we have surjections:
$$H_k({\cal V}-\L) \rightarrow H_k(V_{\cal A}-\L) \ \ \ (k \le n-1) $$
This yields that in the diagram
$$ \matrix {  H_k({\cal V}-V_{\A})  & \longrightarrow  &  H_k(T(B)-V_{\A})
\cr
                     & \searrow &  \downarrow \cr
                     &  & H_k(V_{\cal A}-\L) \cr }$$
the vertical arrow is surjective for $k \le n-1$. Hence the claim
about the eigenvalues of the monodromy follows.

\begin{rem} {\rm
One can restrict the collection of strata $S^j_i$ in theorem
\ref{main} by considering only strata which closures are dense
(cf. \cite{stv}) edges. Recall that a dense edge is an
intersection of hyperplanes of arrangement $L$ such that the
central arrangement $\A_L=\{H \in \A \vert L \subseteq H \}$ is
indecomposable in the sense that the latter arrangement cannot be
split into a union of two subarrangement which can be written in
appropriate coordinates as arrangements of disjoint sets of
variables.

Indeed, a Thom-Sebastiani (cf. \cite{thomseb}) type argument shows
that orders of eigenvalues of decomposable arrangements are least
common factor of the orders of monodromy on each factor. }
\end{rem}

\begin{rem} {\rm One can replace $m^j_i$ by the {\it least common multiple}
of the eigenvalues of the monodromy of Milnor fiber of the arrangement 
which appear in the transversal to the stratum section. Such l.c.m.
is a divisor of the multiplicity $m^j_i$ (cf. section 2) but can be 
smaller than $m^j_i$. The simplest example is two lines through a point.
Here the multiplicity is $2$ but the only  eigenvalues of the monodromy is 
$1$.}
\end{rem}

\begin{rem} {\rm
Theorem \ref{main} gives also restriction on the local systems
corresponding to the homomorphism of the fundamental group sending
each meridian to $1 \over e$ which have non vanishing cohomology
in dimension $k-1$: non vanishing occurs only if $e$ divides at
least one of the multiplicities $m^j_k, j \le k$. }
\end{rem}

\begin{coro} \label{relprime} {\rm Order of the monodromy operator acting on
$H_i(F_A), 1 \le i \le n-1$  for an arrangement of $d$ hyperplanes
in $\PP^n$ divides $d$ and at least one of numbers $m^i_j(H)$ for
each hyperplane $H$ of the arrangement. In particular, if for at
least one hyperplane each $m^j_i(H)$ is relatively prime to $d$,
then eigenvalues different from $1$ appear only in top dimension
(i.e. $n$).}
\end{coro}

\section{Line arrangements}

Let ${\cal A}$ be an arrangement of $d$ lines in ${\bf P}^2$.
This is as a curve of degree $d$
having only ordinary singularities. We shall apply the calculation
of the Alexander module of plane algebraic curves to this
special curve. We refer to \cite{arcata}, \cite{loeser} or \cite{tokyo}
for definition of constants and ideals of quasiadjunction.
Since the singularities of the curve in question are
(weighted) homogeneous one can use
 the description of the constants of quasiadjunction from
\cite{arcata} sect. 5 (cf. also \cite{tokyo} sect.3).
One obtains that in the case of a point of multiplicity $m$ the
constants of quasiadjunction are ${{m-2} \over m}, {{m-3} \over
m},...,{1 \over m}$ and the ideal of germs $\phi$ in the local
ring of the singular point satisfying $\kappa_{\phi} < \alpha$ is
${\cal M}^{m -[\alpha \cdot m]-2}$ where ${\cal M}$ is the maximal
ideal of the singular point.

\begin{theo} \label{lines}
Let $d={\rm Card} \A $ and $m$ be a divisor of $d$.
Let $\sigma_k(m)$ be the superabundance of the curves of degree
$d-3 -{{kd} \over m}$ such that the local equation at a point of
multiplicity $m$ belongs to the ideal ${\cal M}^{m-[{{k \cdot d}
\over m}]-1}$. Then the multiplicity
of an eigenvalue $exp({{2 \pi i k} \over m})$ of the monodromy of
the Milnor fiber acting on $H_1(F_A,{\CC})$ is equal to the
$\sigma_k(m)+\sigma_{d-k}(m)$.

\end{theo}

\par \noindent {\bf Proof.} We shall use the identification
(\ref{nonones}). We have $H^{0,1}(\tilde V_{\A})_{\xi}=
H^{1,0}(\tilde V_{\A})_{\bar \xi}$ and $H^1(\tilde V_{\A})=
H^{0,1}(\tilde V_{\A}) \oplus H^{0,1}(\tilde V_{\A})$.
According to theorem 5.1 in \cite{arcata} we have:
$$H^{0,1}(\tilde V_{\A})_{exp{{2 \pi i k} \over m}}
=H^1({\PP^2},{\cal I}(d-3-{{k\cdot d} \over m}))$$
for the ideal sheaf  ${\cal I}$ defined as follows.
The support of ${\cal O}_{\PP^2}/\cal I$ coincide with
the set of vertices of the arrangement $\A$ and the stalk at
a singular point $P$ consists of germs $\phi \in {\cal O}_P$
such that $\kappa_{\phi} < {k \over m}$.
Now the theorem follows from
the above description of ideals of quasiadjunction of
ordinary singularities.

\section{Remarks and Examples}

\begin{rem} $\zeta$-function of monodromy. {\rm One can easily see the 
the relation:} 
\begin{equation} \label{zeta}
\zeta_{\cal A}(t)=\prod det(1-T_{f,i}t,H_i(V_A))^{(-1)^{i}}
=(1-t^d)^{\chi({\PP^n}-A)}
\end{equation}
{\rm Indeed, we have}  $\zeta_{\cal A}(t)=
\prod det(1-T_{f,i}t,H_i(V_A))^{(-1)^{i}}=
\prod det(1-T_{f,i}t,C_i(V_A))^{(-1)^{i}}$ {\rm where} 
$C_i(V_A)$ {\rm denote the} $i$-{\rm chains. }  
\end{rem}
\begin{exam} Braid arrangement i.e. the arrangement in ${\bf P}^n$
with hyperplanes given by $x_i=x_j (i,j=0,...,n,i \ne j)$ . 
{\rm In the case $n=2$ (three lines in $\PP^2$ passing through a point)
section \ref{lines} yields that the multiplicity of eigenvalue $\omega_3$ 
is equal to $1$. Since $\chi({\bf P}^3-\A)=-1$ and 
${\rm dim} H_1 ({\bf P}^3-\A)=2$ we obtain $(1-t),(1-t)(1-t^3),1$
as the characteristic polynomials of the monodromy acting 
on $H_0,H_1$ and $H_2$ respectively.

In the case
$n=3$ we have the arrangement of $6$ planes with $4$ lines and one
vertex. Each plane contains two strata of codimension 2 in $\PP^3$
(each has a line as a closure) 
having multiplicity 3 and hence the eigenvalues of the monodromy
acting on $H_1$ of the Milnor fiber $\prod (x_i-x_j)=1$ are the
roots of unity of order 3 or 1. Similarly the eigenvalues of the 
monodromy acting on $H_2$ are roots of unity of degree either $1,3$ or $6$
(in fact all orders do occur (cf. \cite{suciu})).

In fact $\omega_3$ {\it is} an eigenvalue of $T_1$. Indeed, 
by Lefschetz type argument
(cf 2.5) eigenvalues are the eigenvalues of the monodromy of the
arrangement of 6 lines in ${\bf P}^2$ formed the lines containing
the sides and the medians of a triangle. One can find the
multiplicity of the eigenvalue of $\omega_3$ as the ${\rm dim}
H^1({\cal I}(1))$ where $\cal I$ is the ideal sheaf of the collection
of triples points in this arrangement of lines. Since these
4 triple points form a complete intersection of two quadrics we
have:
$$0 \rightarrow \O(-4) \rightarrow  \O(-2) \oplus \O(-2)
\rightarrow {\cal I} \rightarrow 0$$
and hence $H^1({\cal I}(1))=H^2({\O(-3)})={\bf C}$.
Therefore the multiplicity of the eigenvalue $\omega_3$
is 1. 

The divisibility theorem has the following consequence. Since a
generic plane section has only triple points the eigenvalues most
have the order 3 or 1 and the order must be a divisor of $
{{n(n-1)} \over 2} $. Hence if $n = 2 \ {\rm mod \ 3} $ the order of
an eigenvalue must be 1.

Similarly, generic section of a braid arrangement by a 3-space has
only points of multiplicity $3$ and $6$. Hence if $n$ is relatively 
prime to $6$ then the eigenvalues of the monodromies $T_1$ and $T_2$
(i.e. on $H_1$ and $H_2$ respectively) are equal to $1$, etc. 

}
\end{exam}

\begin{exam} {\rm More generally, consider an arrangement such 
that only three hyperplanes meet in each edge of codimension two
(cf. previous example). Then if the number
of hyperplanes is not divisible by 3 then the eigenvalues of the
monodromy are equal to 1. In particular the homology of the Milnor
fiber in all dimensions except top coincide with the cohomology of
$\PP^n-\A$.}
\end{exam}

\begin{exam} {\rm Consider a line arrangement having only triple points.
It follows from (\ref{nonones}) that the multiplicities of the
eigenvalues $exp{{2 \pi i } \over 3}$ and  $exp{{4 \pi i } \over 3}$
are the same. It follows from the divisibility theorem
at infinity \cite{duke} that the characteristic polynomial
of the monodromy divides $(t^d-1)^{d-2}(t-1)$. In particular
the multiplicity of the eigenvalue $exp{{2 \pi i } \over 3}$
does not exceed $d-2$ Since the multiplicity of the eigenvalue $1$
is $d-1$ we obtain:
$$H_1(F_A,{\bf C}) =({\bf C}[t,t^{-1}]/(t^3-1))^s \oplus
  ({\bf C}[t,t^{-1}]/(t-1))^{d-1-s}$$
where $s$ is the superabundance of the curves of degree $d-3-{d \over 3}$
passing through the set of vertices of multiplicity $3$. (cf. \cite{dedham})}
\end{exam}

\begin{exam} {\rm Consider an arrangement of $p$ hyperplanes where
$p$ is a prime and such that not all hyperplanes are passing
through a point. Then eigenvalues different from $1$ appear only
in the top dimension. Indeed the multiplicity of any stratum is
less than $p$ and the claim follows from \ref{relprime}}.
\end{exam}

\end{document}